\documentclass[11pt]{article}
\usepackage{amsmath,amsfonts,latexsym,amssymb,graphicx,abstract,here,bbm,url}
\advance \topmargin by -\headheight \advance \topmargin by
-\headsep

\def\l{\lambda}

\def\P{{\mathbb P}}

\def\ai{\mbox{Ai}}

\def\a{\alpha}

\def\R{\mathbb R}

\def\F{{\mathbb F}}

\def\P{{\mathbb P}}

\def\l{\lambda}
\def\labda1{\lambda_1}
\def\labda2{\lambda_2}

\def\e{\varepsilon}
\def\f{\phi}

\def\s{\sigma}

\def\comment#1{\relax}

\def\=in{\mathop{\rm =}}

\def\eop{\hfill\mbox{$\Box$}\newline}

\newtheorem{theorem}{Theorem}[section]

\newtheorem{lemma}{Lemma}[section]

\newtheorem{remark}{Remark}[section]

\begin{document}
\title{Vertices of the least concave majorant of Brownian motion with parabolic drift}
\author{Piet Groeneboom}
\date{\today}
\affiliation{Delft University of Technology,\\
Delft Institute of Applied Mathematics,\\
Mekelweg 4,\\
2628CD Delft,\\
The Netherlands,\\
http://dutiosc.twi.tudelft.nl/\,$_{\widetilde{~}}$pietg/
}
\AMSsubject{60J65}
\maketitle
\begin{abstract}
It was shown in \cite{piet:83} that the least concave majorant of one-sided Brownian motion without drift can be characterized by a jump process with independent increments, which is the inverse of the process of slopes of the least concave majorant. This result can be used to prove the result in \cite{Sparre:54} that the number of vertices of the smallest concave majorant of the empirical distribution function of a sample of size $n$ from the uniform distribution on $[0,1]$ is asymptotically normal, with an asymptotic expectation and variance which are both of order $\log n$.

A similar (Markovian) inverse jump process was introduced in \cite{piet:89}, in an analysis of the least concave majorant of two-sided Brownian motion with a parabolic drift. This process is quite different from the process for one-sided Brownian motion without drift: the number of vertices in a (corresponding slopes) interval has an expectation proportional to the length of the interval and the variance of the number of vertices in such an interval is about half the size of the expectation, if the length of the interval tends to infinity. We prove an asymptotic normality result for the number of vertices in an increasing interval, which translates into a corresponding result for the least concave majorant of an empirical distribution function of a sample of size $n$, generated by a strictly concave distribution function. In this case the number of vertices is of order cube root $n$ and the variance is again about half the size of the asymptotic expectation.

As a side result we obtain some interesting relations between the first moments of the number of vertices, the square of the location of the maximum of Brownian motion minus a parabola, the value of the maximum itself, the squared slope of the least concave majorant at zero, and the value of the least concave majorant at zero. 
\end{abstract}

\newpage
\section{Introduction}
\label{intro}
It was shown in \cite{piet:83} that one-sided Brownian motion can be generated by a jump process with independent increments (which is the inverse of the process of slopes of the least concave majorant) together with Brownian excursions between successive vertices of the least concave majorant. This decomposition of Brownian motion, using the inverse process, was also analyzed in \cite{pitman:83}, where certain path decomposition results, introduced by David Williams, were applied.

The study of the (least) concave majorant of Brownian motion in \cite{piet:83} was actually motivated by the wish to give an alternative derivation of the asymptotic distribution for certain statistics, studied in \cite{behnen:75} and \cite{scholz:83}, and first proved in \cite{piet_ron:83}, where the asymptotic distribution was found by analyzing the spacings, induced by the least concave majorant of the empirical distribution function.

As a side effect, \cite{piet:83} also threw some new light on a result of \cite{Sparre:54}, which is stated below. A straightforward proof of this result, using characteristic functions and the Poisson representation in \cite{piet_ron:83}, is given in \cite{piet_rik:93}.

\begin{theorem}
{\rm [Sparre Andersen (1954)]} Let $N_n$ be the number of vertices of the least concave majorant of the empirical distribution function of a sample of size $n$ from the uniform distribution on $[0,1]$. Then
$$
\frac{N_n-\log n}{\sqrt{\log n}}\stackrel{{\cal D}}\longrightarrow N(0,1),
$$
where $\stackrel{{\cal D}}\longrightarrow$ denotes convergence in distribution, and $N(0,1)$ is the standard normal distribution.
\end{theorem}

The corresponding result for the Brownian bridge on $[0,1]$, which follows from \cite{piet:83}, is:

\begin{theorem}
Let $C$ be the least concave majorant of the Brownian bridge on $[0,1]$. Then the number of vertices $N_n$ of $C$ on the interval $[1/n,1-1/n]$ satisfies:
$$
\frac{N_n-\log n}{\sqrt{\log n}}\stackrel{{\cal D}}\longrightarrow N(0,1).
$$
\end{theorem}

If one studies more closely ``where the action is", in the sense that the number of vertices increases to infinity, it turns out that all the action is near $0$ and $1$: on an interval $[\e,1-\e]$, where $\e>0$, there will, with probability one, only be finitely many vertices.

The situation is strikingly different for the least concave majorant of two-sided Brownian motion minus a parabola. Here the action is ``the same everywhere", and the point process of locations of vertices is stationary. The real purpose of the papers \cite{piet:85} and \cite{piet:89} was to analyze the (stationary) process
\begin{equation}
\label{V_process}
\left\{V(a)-a:a\in\R\right\},
\end{equation}
where
$$
V(a)=\mbox{argmax}_{t\in\R}\left\{W(t)-(t-a)^2\right\},\,a\in\R,
$$
and $W$ is standard two-sided Brownian motion, originating from zero. The process $V$ itself is a pure jump process, which runs through the locations of the vertices of the least concave majorant.
At points where $V(a)$ is not uniquely determined (which happens if $a$ is a slope of the least concave majorant), we take the largest value $t$ at which $W(t)-(t-a)^2$ is maximal. In this way the process $V$ becomes right-continuous.

\begin{figure}[!ht]
\begin{center}
\includegraphics[scale=0.5]{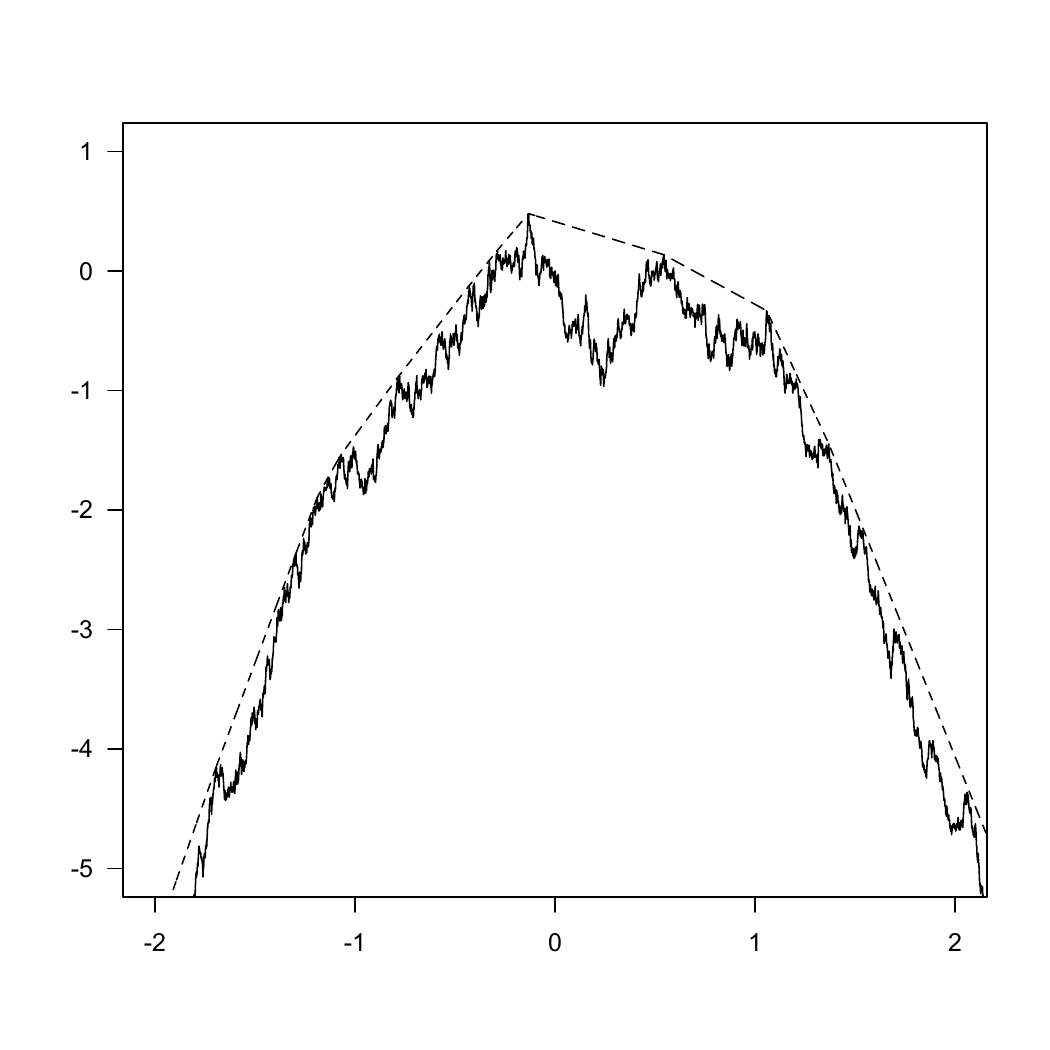}
\end{center}
\caption{The least concave majorant of $W(t)-t^2$.}
\label{fig:conc_maj}
\end{figure}

The infinitesimal generator of the process (\ref{V_process}) is given in Theorem 4.1 of \cite{piet:89}, where it is expressed in terms of Airy functions. However, most attention has been for the result on the distribution of $V(0)$, which gave an analytic expression for the limit distribution of a whole class of so-called ``isotonic estimators", for example the pointwise limit distribution of an estimator of the mode, discussed in \cite{chernoff:64}, and the pointwise limit behavior of the Grenander (maximum likelihood) estimator of a decreasing density, see, e.g., \cite{prakasa:69} and \cite{piet:85}.

For results on global functionals, however, like the $L_1$ or $L_2$ distance of the Grenander estimator to the underlying density, or the number of its jumps, one needs information on the whole process $V$, and not only on its pointwise behavior. In this paper we will show how one can extract information from Theorem 4.1 in \cite{piet:89} in the derivation of a central limit theorem for the number of points jump of $V$ in an increasing interval. The result has a rather large number of applications in statistics, but we will only sketch one such result for the number of points of jump of the Grenander estimator (which is equivalent to the corresponding result for the number of vertices of the least concave majorant).

Our main result is the following central limit theorem, which is proved at the end of section \ref{section:BM}.

\begin{theorem}
\label{th:asymp_normality_BM}
Let $N[a,b]$ be the number of jumps of the process $V$ in the interval $[a,b]$. Then
$$
\frac{N[a,b]-k_1(b-a)}{\sqrt{k_2(b-a)}}\stackrel{{\cal D}}\longrightarrow N(0,1),\mbox{ as }b-a\to\infty,
$$
where $k_1\approx2.10848$ and $k_2\approx1.029$, and $N(0,1)$ is the standard normal distribution.
\end{theorem}

Perhaps somewhat remarkably, the difference between the results for the least concave majorants of one-sided Brownian motion without drift and two-sided Brownian motion with a parabolic drift has its counterpart in the difference between the convex hulls of uniform samples of points from the interior a convex polygon and from the interior of a convex figure with a smooth boundary, see \cite{piet:88}. In this case one also meets the rates $\log n$ and $n^{1/3}$ for the number of vertices of the convex hulls of the samples, with corresponding central limit results.

\section{The number of jumps of the process $V$ in an increasing interval}
\label{section:BM}
 Although \cite{piet:89} has the simpler conceptual characterization, the characterization in \cite{piet:85} might be more useful for  numerical computations. It was also used in \cite{piet_jon:01}, where Chernoff's density and its moments were computed. Chernoff's density is the density of $V(0)$, which often occurs as limit of isotonic estimators and in particular in the limit distribution of the Grenander estimator. Here, however, we  take \cite{piet:89} as our starting point.

The process (\ref{V_process}) is completely characterized by Theorem 4.1 of \cite{piet:89}. As a corollary we have the following result for the jump measure.

\begin{theorem}
\label{th:jump_measure}
We have, if $y>x$,
$$
\lim_{h\downarrow0} h^{-1}\P\left\{V(a+h)\in a+dy|V(a)=a+x\right\}
=\frac{2(y-x)g(y)p(y-x)}{g(x)}\,dy,
$$
where $g$ has Fourier transform
\begin{equation}
\label{g_Fourier}
\hat g(u)=\int e^{iux} g(x)\,dx=\frac{2^{1/3}}{\mbox{\rm Ai}\left(i2^{-1/3}u\right)}\,,\,u\in\R,
\end{equation}
and
$$
p(u)=p_0(u)+\left(2\pi u^3\right)^{-1/2},\,u>0,
$$
where $p_0$ has Laplace transform
$$
\hat p_0(u)=\frac{2^{2/3}\mbox{\rm Ai}\,'\bigl(2^{-1/3}u\bigr)}{\mbox{\rm Ai}\bigl(2^{-1/3}u\bigr)}+2^{1/2}\sqrt{u}.
$$
\end{theorem}

\begin{remark}
{\rm Note that we define the Fourier transform in the ``probabilistic way", in analogy with the definition of the characteristic function of a probability distribution.
}
\end{remark}

\begin{figure}[!ht]
\begin{center}
\includegraphics{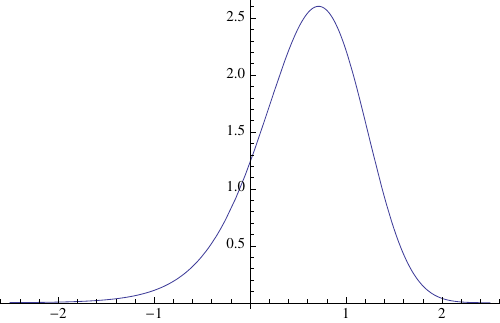}
\end{center}
\caption{The function $g$.}
\label{fig:g}
\end{figure}

A picture of the function $g$, using the representation
$$
g(x)=\frac1{2^{2/3}\pi}\int_{-\infty}^{\infty}\frac{e^{-i u x}}{\mbox{\rm Ai}(i2^{-1/3}u)}\,du,
$$
which follows from (\ref{g_Fourier}), is shown in Figure \ref{fig:g}.

\begin{figure}[!ht]
\begin{center}
\includegraphics{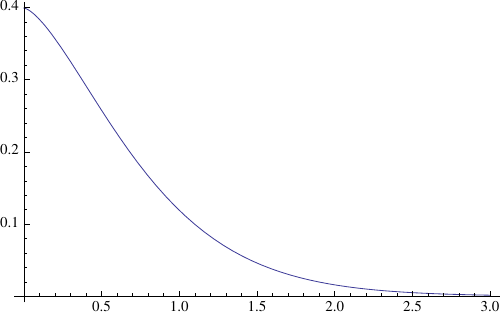}
\end{center}
\caption{The function $u\mapsto u^{3/2}p(u),\,u>0$.}
\label{fig:p}
\end{figure}

The function $p$ has the representation
\begin{equation}
\label{p_representation}
p(u)=2\sum_{n=1}^{\infty} e^{2^{1/3} \tilde a_n u},\,u>0,
\end{equation}
where the $\tilde a_n$ are the zeros of the Airy function $\ai$ on the negative halfline, see (4.12) in \cite{piet:89}. This expansion is divergent at zero, however, where we have:
\begin{equation}
\label{p_atzero}
p(u)\sim(2\pi u^3)^{-1/2},\,u\downarrow0,
\end{equation}
see part (ii) of Lemma 4.2 in \cite{piet:89}.
This is the reason for considering the regularization
$$
p_0(u)=p(u)-(2\pi u^3)^{-1/2}\,,
$$
and for only using the representation (\ref{p_representation}) for $u\ge1$. If $u<1$ we use the representation given below in (\ref{tilde_p}) of Lemma \ref{lemma:2_p's}.
The function $u\mapsto u^{3/2}p(u)$ is shown in Figure \ref{fig:p}.

For later purposes we summarize in the following lemma some properties of functions $p$ and $g$ and the random variable $V(0)$. This also gives still another regularization of the function $p$.

\begin{lemma}
\label{lemma:p_regularization}
Let the function $p$ be defined as in Theorem \ref{th:jump_measure}. Then:
\begin{enumerate}
\item[(i)]
The function $p_1$, defined by
\begin{equation}
\label{eq:p_regularization}
p_1(x)=\left\{\begin{array}{lll}
xp(x)\,&,\,x>0,\\
0 &,\,x<0,
\end{array}
\right.
\end{equation}
has Fourier transform
\begin{equation}
\label{eq:fourier_transf_p1}
\hat p_1(u)=iu+\frac{2^{1/3}\mbox{\rm Ai}\,'\bigl(-2^{-1/3}iu\bigr)^2}{\mbox{\rm Ai}\bigl(-2^{-1/3}iu\bigr)^2}.
\end{equation}
\item[(ii)] The function $h$, defined by
\begin{equation}
\label{eq:h}
h(x)=\int_{u=0}^{\infty}g(x+u)up(u)\,du,\,x\in\R,
\end{equation}
has Fourier transform
\begin{equation}
\label{eq:fourier_transf_h}
\hat h(u)=-\frac{2^{1/3}iu}{\mbox{\rm Ai}\bigl(2^{-1/3}iu\bigr)}+\frac{2^{2/3}\mbox{\rm Ai}\,'\bigl(2^{-1/3}iu\bigr)^2}{\mbox{\rm Ai}\bigl(2^{-1/3}iu\bigr)^3}.
\end{equation}
\item[(iii)] The random variable $V(0)$ has characteristic function
\begin{align}
E e^{itV(0)}=\frac{1}{2\pi}\int_{u=-\infty}^{\infty}\frac{du}{\mbox{\rm Ai}\left(iu\right)\mbox{\rm Ai}\left(i(2^{-1/3}t+u)\right)},\,t\in\R.
\end{align}
\item[(iv)] The random variable $V(0)$ has expectation zero and second moment
\begin{equation}
\label{eq:2nd_moment_V}
EV(0)^2=-\frac{2^{-2/3}}{6\pi}\int_{u=-\infty}^{\infty}\frac{iu}{\mbox{\rm Ai}\left(iu\right)^2}\,du\approx0.26355964.
\end{equation}
\end{enumerate}
\end{lemma}

\noindent
{\bf Proof.} Part (i): this immediately follows from Theorem \ref{th:jump_measure}, noting that the function $x\mapsto 1/\sqrt{2\pi x}$ has Laplace transform $1/\sqrt{2u}$ and by switching from Laplace transform to Fourier transform.\\
Part (ii): this follows from (i) and Theorem \ref{th:jump_measure}, by noting that the convolution turns into the product of the Fourier transforms (with an added change of sign for the Fourier transform of $p_1$).\\
Part (iii): this is the Fourier transform of the function
$$
x\mapsto \tfrac12g(-x)g(x),\,x\in\R,
$$
which turns into the convolution of the Fourier transform of $x\mapsto g(-x)$ and Fourier transform of $x\mapsto g(x)$ times $(4\pi)^{-1}$.\\
Part (iv) follows from the formula
\begin{align*}
EV(0)^2&=-\frac{2^{-2/3}}{2\pi}\frac{d^2}{dt^2}\int_{u=-\infty}^{\infty}\frac{du}{\mbox{\rm Ai}\left(iu\right)\mbox{\rm Ai}\left(i(t+u)\right)}\Biggr|_{t=0}
=-\frac{2^{-2/3}}{2\pi}\int_{u=-\infty}^{\infty}\left\{\frac{iu}{\ai(iu)^2}-\frac{2\ai\,'(iu)^2}{\ai(iu)^4}\right\}\,du\\
&=-\frac{2^{-2/3}}{2\pi}\int_{u=-\infty}^{\infty}\left\{\frac{iu}{\ai(iu)^2}-\frac{2iu}{3\ai(iu)^2}\right\}\,du
=-\frac{2^{-2/3}}{6\pi}\int_{u=-\infty}^{\infty}\frac{iu}{\ai(iu)^2}\,du,
\end{align*}
where we use integration by parts and the Airy equation for the term involving $\mbox{Ai}'$.
 \eop

\begin{remark}
\label{remark_on_p}
{\rm
The function $p$ does not have the same meaning in \cite{piet:85}. If we denote the function $p$ of \cite{piet:85} by $\tilde p$, we have:
$$
p(t)=\frac{e^{t^3/6}}{\sqrt{2\pi}}\left\{\tilde p(t)+t^{-3/2}\right\},\,t>0.
$$
The function $\tilde p$ has an expansion which converges at zero.
}
\end{remark}

\begin{lemma}
\label{lemma:2_p's}
The function $p$ can be written
$$
p(t)=\frac{e^{t^3/6}}{\sqrt{2\pi}}\left\{\tilde p(t)+t^{-3/2}\right\},
$$
where
\begin{eqnarray}
\label{tilde_p}
\tilde p(t) = 
 -\sqrt{\frac{\pi}{2}} \sum_{k=0}^{\infty} a_k t^{3k} 
  + \sum_{k=1}^{\infty} b_k t^{3(k-1/2)}\,,\,t>0,
\end{eqnarray} 
and the coefficients $a_k$ and $b_k$ are recursively defined as follows.
Set $c_0 =1$ and 
$$
c_n = - 2^{-4} \frac{(2n-3)(2n+1)}{n^2(2n-1)} c_{n-1} \, , \qquad 
n = 1, 2, \ldots \, .
$$
Then with $a_0 =1$, $b_1 = 2/3$, and 
$B(p,q) \equiv \Gamma(p) \Gamma (q)/ \Gamma (p+q)$, the 
standard Beta function, set
\begin{equation}
\label{a-coefficients}
a_n = c_n - \sum_{k=0}^{n-1} \frac{1}{\pi k!(-2)^k} b_{n-k} 
B(3n-2k-1/2 , k+3/2) \, , \qquad 
n = 1, 2, \ldots \, ;
\end{equation}
\begin{equation}
\label{b-coefficients}
b_n = \sum_{k=0}^{n-1} \frac{1}{k!(-2)^{k+1}} a_{n-k-1} 
B(3n-2k-2 , k+3/2) \, , \qquad 
n = 2, 3, \ldots \, .
\end{equation}
\end{lemma}

\noindent
{\bf Proof.} This follows from Remark \ref{remark_on_p} and Theorem 4.2 in \cite{piet:85}.\eop

Let $u_2: \R \rightarrow \R$ be defined by 
\begin{eqnarray}
\label{function_u2_1}
u_2(x) 
& = & 2x - \frac{1}{\sqrt{2 \pi}} \int_0^{\infty} \tilde p(y) \
      \exp(-\tfrac12 y(2x+y)^2 ) dy \nonumber\\
& & \quad + \frac1{\sqrt{2 \pi}}
    \int_0^{\infty}\left\{ 4x^2+8xy+3y^2\right\} 
    \exp(- \tfrac12 y (2x+y)^2 ) \,y^{-1/2}\, dy 
\end{eqnarray}
if $x \in [-1, \infty)$, and 
\begin{equation}
\label{function_u2_2}
u_2(x) = \exp\left(\tfrac23 x^3\right) 4^{1/3} \sum_{k=1}^{\infty} 
       \exp(-2^{1/3} \tilde{a}_k x )/\ai'(\tilde{a}_k) 
\end{equation}
if $x \in (-\infty , -1]$; here $\ai'$ is the derivative of the Airy
function $\ai$.

The notation $u_2$ is used because $u_2$ has the interpretation
$$
u_2(t)=\lim_{x\uparrow t^2}\frac{\partial}{\partial x}u(t,x),
$$
where $u(t,x)$  is the solution of the heat equation
$$
\frac{\partial}{\partial t}u(t,x)=-\tfrac12\frac{\partial^2}{\partial x^2}u(t,x),
$$
for $x\le t^2$, under the boundary conditions
$$
u(t,t^2)\stackrel{\mbox{def}}=\lim_{x\uparrow t^2}u(t,x)=1,\qquad \lim_{x\downarrow -\infty}u(t,x)=0,\qquad t\in\R.
$$
This function occurred in the paper \cite{chernoff:64}, where the density of the location of the maximum of Brownian motion minus a parabola was first characterized.

The following result summarizes the correspondence between the results in \cite{piet:85} and \cite{piet:89}.

\begin{theorem}
\label{th:jump_measure1}
\begin{enumerate}
\item[(i)]
We have, if $y>x$,
\begin{align*}
&\lim_{h\downarrow0} h^{-1}\P\left\{V(a+h)\in a+dy|V(a)=a+x\right\}\\
&=\frac{2(y-x)u_2(y)e^{-\tfrac12(y-x)(x+y)^2}\left\{\tilde p(y-x)+(y-x)^{-3/2}\right\}}{u_2(x)\sqrt{2\pi}}\,dy,
\end{align*}
where $u_2$ is defined by (\ref{function_u2_1}) and (\ref{function_u2_2}) and $\tilde p$ by (\ref{tilde_p}).
\item[(ii)]
$$
u_2(x)=e^{\tfrac23x^3}g(x),\,x\in\R\quad\mbox{ and }\quad\tilde p(t)=\sqrt{2\pi}e^{-t^3/6}p(t)-t^{-3/2},\,t>0,
$$
where the functions $g$ and $p$ are defined as in Theorem \ref{th:jump_measure}.
\item[(iii)] The density of $V(0)$ is given by
$$
f_{V(0)}(x)=\tfrac12u_2(x)u_2(-x)=\tfrac12g(x)g(-x),\,x\in\R,
$$
where $g$ is defined as in Theorem \ref{th:jump_measure}.
\end{enumerate}
\end{theorem}

\begin{remark}
\label{remark_on_piet_jon}
{\rm
In \cite{piet_jon:01} part (iii) of Theorem \ref{th:jump_measure1} was used in the computation of the Chernoff distribution. The function $u_2$ corresponds to the function $k_1$ in \cite{piet:89} and part (i) of of Theorem \ref{th:jump_measure1} corresponds to the first version of the infinitesimal generator of the process $V(a)-a$, given in Theorem 4.1 of that paper.

It is seen from Theorem \ref{th:jump_measure1} and (\ref{function_u2_1}) that the function $p$ (or alternatively, the regularization $\tilde p$) is the fundamental function; both the jump measure and the  density of $V(0)$ are expressed in terms of $p$.
}
\end{remark}

We are now ready to compute the expectation of the number of jumps of the process $V$ in an interval (of slopes) $[a,b]$, where we use similar techniques as in \cite{piet:88}, which dealt with convex hulls of samples of points from the interior of a convex set in the plane.

Let the function $\f$ be defined by
\begin{equation}
\label{phi_b}
\f(x)=2\int_0^{\infty} \frac{g(x+u) u p(u)}{g(x)}\,du,\,x\in\R,
\end{equation}
$\f(x)$ is the integrated jump measure, starting from position $x$.
Moreover, let $N[a,b]$ denote the number of jumps of the process $V$ in the interval $[a,b]$.
Then Theorem \ref{th:jump_measure} tells us that
\begin{equation}
\label{martingale}
b\mapsto N[a,b]-\int_a^b \f(V(c)-c)\,dc,\,b\ge a,
\end{equation}
is a martingale w.r.t.\ the filtration, generated by $V(b),\,b\ge a$. As a consequence, we have the following result.

\begin{lemma}
\label{lemma:expectation}
Let $N[a,b]$ be the number of jumps of the process $V$ in the interval $[a,b]$. Then
$$
E N[a,b]=k_1(b-a),
$$
where
\begin{equation}
\label{constant_k1}
k_1=\int_{-\infty}^{\infty} g(-x)\,dx\int_{y=x}^{\infty}g(y)(y-x)p(y-x)\,dy\approx2.10848.
\end{equation}
\end{lemma}

\noindent
{\bf Proof.} We get from Theorem \ref{th:jump_measure}:
$$
E N[a,b]=\int_a^b E\f(V(c)-c)\,dc.
$$
Using the stationarity of the process $c\mapsto V(c)-c$ we get:
\begin{align*}
&\int_a^b E\f(V(c)-c)\,dc
=\int_a^b E\f(V(0))\,dc\\
&=(b-a)\int_{-\infty}^{\infty}f_{V(0)}(x)\,dx\int_{x}^{\infty} \frac{2g(y)(y-x)p(y-x)}{g(x)}\,dy\\
&=(b-a)\int_{-\infty}^{\infty}\tfrac12g(x)g(-x)\,dx\int_{x}^{\infty} \frac{2(y-x)g(y)p(y-x)}{g(x)}\,dy\\
&=(b-a)\int_{-\infty}^{\infty}g(-x)\,dx\int_{x}^{\infty} g(y)(y-x) p(y-x)\,dy,
\end{align*}
and the result follows. The constant $k_1$ was determined numerically by using Theorem \ref{th:jump_measure1}.\eop

\begin{remark}
{\rm As one of the referees remarks, Fourier analysis, applied on the right-hand side of (\ref{constant_k1}), gives:
\begin{align}
\label{2nd_representation_k1}
k_1&=-\frac{2^{5/3}}{6\pi}\int_{u=-\infty}^{\infty}\frac{iu}{\mbox{Ai}\bigl(2^{-1/3}iu\bigr)^2}\,du=-\frac{2^{7/3}}{6\pi}\int_{u=-\infty}^{\infty}\frac{iu}{\mbox{Ai}\bigl(iu\bigr)^2}\,du\nonumber\\
&=8E V(0)^2=\tfrac83E\max_{t\in\R}\left(W(t)-t^2\right),
\end{align}
where $W$ is standard two-sided Brownian motion, originating from zero, and $V(0)$ is defined as in (\ref{V_process}), for $a=0$. This follows from Lemma \ref{lemma:p_regularization}, since we get by Parseval's formula:
\begin{align*}
&\int_{-\infty}^{\infty} g(-x)\,dx\int_{y=x}^{\infty}g(y)(y-x)p(y-x)\,dy
=\frac1{2\pi}\int_{u=-\infty}^{\infty} \overline{\hat g(-u)}\hat h(u)\,du\\
&=\frac1{2\pi}\int_{u=-\infty}^{\infty}\left\{\frac{2\mbox{Ai}'\bigl(2^{-1/3}iu\bigr)^2}{\mbox{Ai}\bigl(2^{-1/3}iu\bigr)^4}-\frac{i2^{2/3}u}{\mbox{Ai}\bigl(2^{-1/3}iu\bigr)^2}\right\}\,du=-\frac{2^{5/3}}{6\pi}\int_{u=-\infty}^{\infty}\frac{iu}{\mbox{Ai}\bigl(2^{-1/3}iu\bigr)^2}\,du,
\end{align*}
where we use integration by parts and the Airy equation for the term involving $\mbox{Ai}'$ as in the proof of part (iv) of Lemma \ref{lemma:p_regularization}.

This also gives an interesting relation between the moments of the location of the maximum and moments of the maximum itself. By \cite{janson:10} we get:
$$
E\max_{t\in\R}\left(W(t)-t^2\right)=2^{-1/3}EM=0.790679,
$$
where $M$ is the maximum of $W(t)-\tfrac12t^2$, see (1.7) and (2.5) of their paper, which is in accordance with the value, given in Lemma \ref{lemma:expectation}. The integral representation for the maximum of $W(t)-t^2$ of type (\ref{2nd_representation_k1}) above corresponds to (2.1) in their paper (after replacing $W(t)-t^2$ by $W(t)-\tfrac12t^2$, see also Remark \ref{remark:scaling_constants} below). The value $EV(0)^2$ was computed in \cite{piet_jon:01}, where it is given by $0.26355964$ (note that this is also given in part (iv) of Lemma \ref{lemma:p_regularization}), and this gives $k_1=2.10848$ again. For convenience, we state this in a separate lemma.
}
\end{remark}

\vspace{0.3cm}
\begin{lemma}
The constant $k_1$ in Lemma \ref{lemma:expectation} has the representation
\begin{equation}
\label{eq:k_1_representation}
k_1=8 EV(0)^2=\tfrac83E\max_{t\in\R}\left(W(t)-t^2\right)\approx2.10848.
\end{equation}
\end{lemma}

\vspace{0.3cm}
\begin{remark}
\label{remark:scaling_constants}
{\rm
Note that
\begin{align*}
&\mbox{argmax}_{t\in\R}\left\{W(t)-ct^2\right\}
\stackrel{{\cal D}}=\mbox{argmax}_{t\in\R}\left(c^{-1/3}\left\{W(c^{2/3}t)-\left(c^{2/3}t\right)^2\right\}\right)\\
&=\mbox{argmax}_{t\in\R}\left\{W(c^{2/3}t)-\left(c^{2/3}t\right)^2\right\}
=c^{-2/3}\mbox{argmax}_{u\in\R}\left\{W(u)-u^2\right\},
\end{align*}
which implies that, if we define $k_1(1)=k_1$ and $k_2(1)=k_2$, and denote the corresponding constants for the process $t\mapsto W(t)-ct^2$ by $k_1(c)$ and $k_2(c)$:
$$
k_i(c)=c^{2/3}k_i,\,i=1,2.
$$
Relation (\ref{2nd_representation_k1}) changes into:
\begin{equation}
\label{3rd_representation_k1}
k_1(c)=8c^2EV_c(0)^2=\frac{8c}3E\max_{t\in\R}\left(W(t)-ct^2\right),
\end{equation}
where $a\mapsto V_c(a)$ is the process of locations of maxima of $t\mapsto W(t)-c(t-a)^2$.
}
\end{remark}

\begin{remark}
\label{remark:Woodroofe}
{\rm As also pointed out by one of the referees, \cite{meyer_woodroofe:00} represent the constant $k_1(1/2)$,
where $k_1(1/2)$ is defined as in Remark \ref{remark:scaling_constants}, in their Corollary 4 as the sum of two constants:
\begin{equation}
\label{meyer_woodroofe}
k_1(1/2)=E\tilde X(0)+E\tilde X'(0)^2,
\end{equation}
where $X(t)=W(t)-\tfrac12t^2$ and $\tilde X$ is the greatest concave majorant of $X$, with slope $\tilde X'(0)$ at zero. (I switch from the convex minorants of Brownian motion plus a parabola to the concave majorants of Brownian motion minus this parabola here; this gives the same $k_1$). \cite{meyer_woodroofe:00} give simulation results for $X(t)=W(t)-\tfrac12t^2$, which would imply that $k_1(1/2)\approx1.289$.
We get from Remark \ref{remark:scaling_constants}:
$$
k_1(1/2)=2^{-2/3}k_1(1)\approx1.32826,
$$
which is larger than the value arising out of the simulations in \cite{meyer_woodroofe:00}, indicating that it is very hard to obtain precise values of these constants by direct simulation of Brownian motion.

Denoting the least concave majorant of the process $t\mapsto W(t)-ct^2$ by $\tilde X_c$, and using a notation similar to the notation of Remark \ref{remark:scaling_constants}, we would get the relation
\begin{equation}
\label{meyer_woodroofe2}
k_1(c)=2cE\tilde X_c(0)+E\tilde X_c'(0)^2,
\end{equation}
which indeed is compatible with the relation:
$$
k_1(c)=c^{2/3}k_1(1).
$$
}
\end{remark}

\begin{remark}
\label{remark:Woodroofe2}
{\rm There exists a simple relation between $V_c(0)$ and $\tilde X'_c(0)$, where $V_c$ is defined as in Remark \ref{remark:scaling_constants} and $\tilde X'_c(0)$ is defined as in Remark \ref{remark:Woodroofe}. This follows from the so-called ``switch relation":
$$
\tilde X'_c(0)\le 2ca\iff V_c(a)\ge 0.
$$
Since
$$
\P\left\{V_c(a)\ge 0\right\}=\P\left\{V_c(a)-a\ge -a\right\}=\P\left\{V_c(0)\ge -a\right\}=\P\left\{2cV_c(0)\le 2ca\right\},
$$
we get that $\tilde X'_c(0)$ and $2cV_c(0)$ have the same distribution, and hence:
$$
E\tilde X'_c(0)^2=4c^2EV_c(0)^2.
$$
Combining the preceding remarks, and in particular assuming that (\ref{meyer_woodroofe2}) (or (\ref{meyer_woodroofe})) holds, we obtain:
$$
k_1(c)=2cE\tilde X_c(0)+4c^2EV_c(0)^2=8c^2EV_c(0)^2,
$$
which implies
$$
E\tilde X_c(0)=2cEV_c(0)^2=\tfrac23E\max_{t\in\R}\left\{W(t)-ct^2\right\}.
$$
}
\end{remark}

\vspace{0.3cm}
The asymptotic behavior of the functions $p$ and $\f$ is given in the following lemma.

\begin{lemma}
\label{lemma:p}
We have:
\begin{enumerate}
\item[(i)]
$$
p(t)\sim 2e^{2^{1/3}\tilde a_1t},\,t\to\infty,
$$
where $\tilde a_1$ is the first zero of the Airy function $\ai$ on the negative halfline.
\item[(ii)]
$$
p(t)\sim (2\pi t^3)^{-1/2},\,t\downarrow0.
$$
\item[(iii)]
\begin{align*}
\f(t)&\sim 2t^2,\,t\to-\infty,\qquad \f(t)\sim\frac1{t},\,t\to\infty.
\end{align*}
\end{enumerate}
\end{lemma}

\vspace{0.3cm}
\noindent
{\bf Proof.}
(i) follows from (\ref{p_representation}) and (ii) is the same as (\ref{p_atzero}), which follows from (4.17) of Theorem 4.1 in \cite{piet:89}.\\
(iii). The function $g$ is denoted by $g_1$ in \cite{piet:89}, and hence, according to part (i) of Corollary 3.4, \cite{piet:89}:
$$
g(t)=4^{1/3}\sum_{n=1}^{\infty}\frac{e^{-2^{1/3}\tilde a_n|t|}}{\ai'(\tilde a_n)}\,,\,t<0,
$$
where the $\tilde a_n$ are the zeros of the Airy function on the negative halfline. 

We now have, using part (i) of Corollary 3.4 and part (ii) of Lemma 4.2, \cite{piet:89}, if $t<0$,
\begin{align*}
&\int_0^{|t|}\frac{2up(u)g(t+u)}{g(t)}\,du
=
4\int_0^{|t|}\frac{\sum_{n=1}^{\infty}e^{2^{1/3}\tilde a_n |t+u|}/{\rm Ai}'(\tilde a_n)}{\sum_{n=1}^{\infty}e^{2^{1/3}\tilde a_n |t|}/{\rm Ai}'(\tilde a_n)}\sum_{n=1}^{\infty}e^{2^{1/3}\tilde a_nu}u\,du\\
&\sim 4\int_0^{|t|}\frac{\sum_{n=1}^{\infty}e^{2^{1/3}\{\tilde a_1-\tilde a_n\}(t+u)}/{\rm Ai}'(\tilde a_n)}{\sum_{n=1}^{\infty}e^{2^{1/3}\{\tilde a_1-\tilde a_n\}t}/{\rm Ai}'(\tilde a_n)}\sum_{n=1}^{\infty}e^{-2^{1/3}\{\tilde a_1-\tilde a_n\}u}u\,du\\
&\sim 4{\rm Ai}'(\tilde a_1)\int_0^{|t|}\sum_{n=1}^{\infty}\frac{e^{2^{1/3}\{\tilde a_1-\tilde a_n\}(t+u)}}{{\rm Ai}'(\tilde a_n)}\sum_{n=1}^{\infty}e^{-2^{1/3}\{\tilde a_1-\tilde a_n\}u}u\,du\\
&=4\int_0^{|t|}\left\{1+\sum_{n=2}^{\infty}e^{2^{1/3}\{\tilde a_1-\tilde a_n\}(t+u)}\frac{{\rm Ai}'(\tilde a_1)}{{\rm Ai}'(\tilde a_n)}\right\}
\left\{1+\sum_{n=2}^{\infty}e^{-2^{1/3}\{\tilde a_1-\tilde a_n\}u}\right\}u\,du
\end{align*}
\begin{align*}
&\sim4\int_0^{|t|}
u\left\{1+\sum_{n=2}^{\infty}e^{-2^{1/3}\{\tilde a_1-\tilde a_n\}u}\right\}\,du
+4\int_0^{|t|}
u\sum_{n=2}^{\infty}e^{2^{1/3}\{\tilde a_1-\tilde a_n\}(t+u)}\frac{{\rm Ai}'(\tilde a_1)}{{\rm Ai}'(\tilde a_n)}\,du\\
&\sim4\int_0^{|t|}
u\,du= 2t^2,\,t\to-\infty.
\end{align*}
Furthermore,
\begin{align*}
&\int_{|t|}^{\infty}\frac{2up(u)g(t+u)}{g(t)}\,du=\int_0^{\infty}\frac{2(|t|+u)p(|t|+u)g(u)}{g(t)}\,du\\
&\sim 4^{2/3}\ai'(\tilde a_1)\int_0^{\infty}e^{-2^{1/3}\tilde a_1|t|}\sum_{n=1}^{\infty}e^{2^{1/3}\tilde a_n\{|t|+u\}}\left\{|t|+u\right\} g(u)\,du
=O(|t|),\,t\to-\infty,
\end{align*}
and the first part of (iii) now follows.\\
Using Laplace's method and part (ii) of Corollary 3.4, \cite{piet:89}, we find:
\begin{align*}
&\int_0^{\infty}\frac{2up(u)g(t+u)}{g(t)}\,du
\sim 2\int_0^{\infty}\frac{up(u)(t+u)\exp\left\{-\tfrac23(t+u)^3\right\}}{t\exp\left\{-\tfrac23t^3\right\}}\,du\\
&\sim\frac2{\sqrt{2\pi}}\int_0^{\infty}\frac{u^{-1/2}(t+u)\exp\left\{-2t^2u-2tu^2-\tfrac23u^3\right\}}{t}\,du\\
&\sim\frac2{\sqrt{2\pi}}\int_0^{\infty}\frac{u^{-1/2}(t+u)\exp\left\{-2t^2u\right\}}{t}\,du
\sim\frac1{t}\,,\,t\to\infty.
\end{align*}
\eop

\begin{figure}[!ht]
\begin{center}
\includegraphics[scale=0.5]{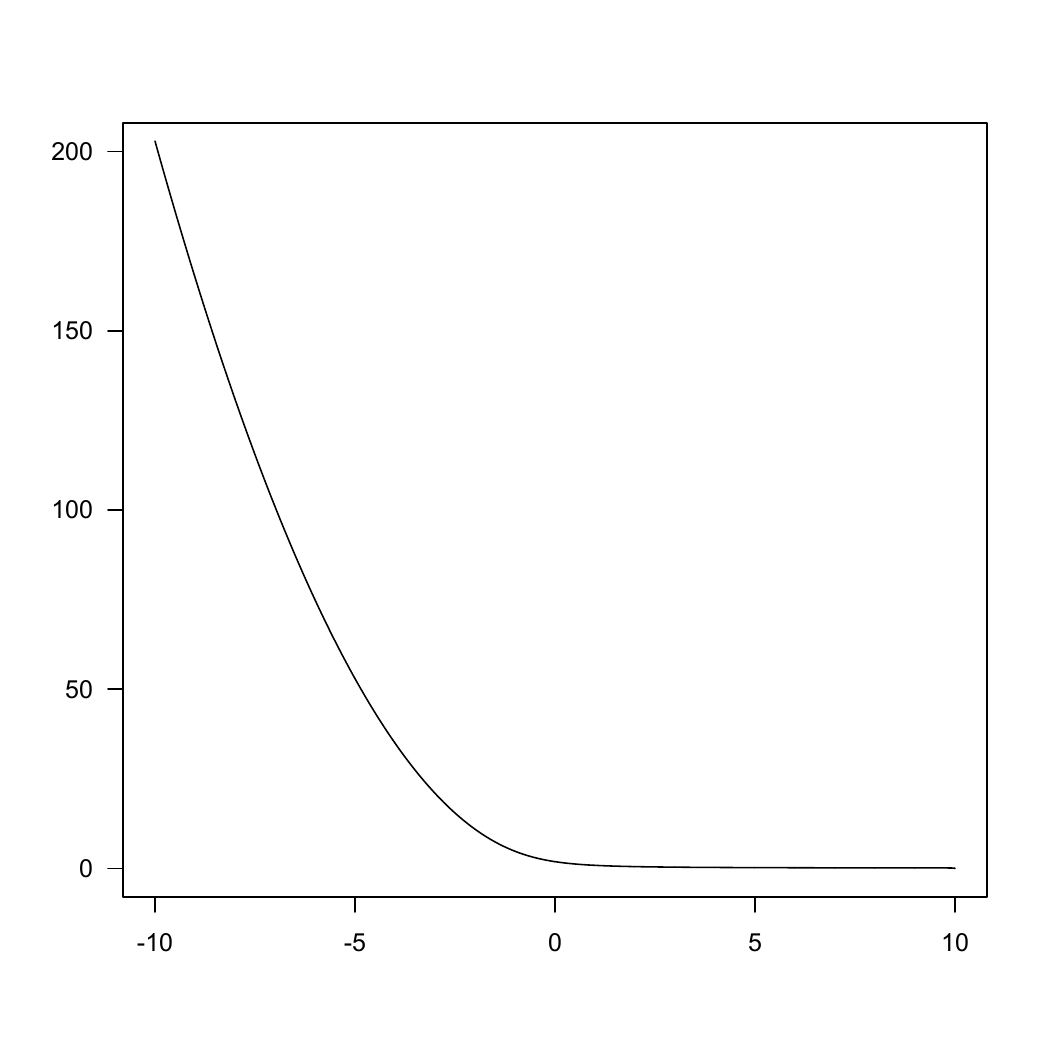}
\end{center}
\caption{The function $\phi$.}
\label{fig:phi}
\end{figure}

\vspace{0.3cm}
For the variance of $N[a,b]$ we get the following result.

\begin{theorem}
\label{th:variance}
Let $N[a,b]$ be the number of jumps of the process $V$ in the interval $[a,b]$. Then:
\begin{equation}
\label{variance_N}
\mbox{\rm var}\bigl(N[a,b]\bigr)
=EN[a,b]+2\int_{a<c_1<c_2<b}\mbox{\rm covar}\left(\f(-V\left(c_1)+c_1\right),\f\left(V(c_2)-c_2\right)\right) dc_1dc_2.
\end{equation}
\end{theorem}

\noindent
{\bf Proof.} We have:
\begin{align*}
&EN(a,b+h)^2-EN[a,b]^2
=E N(b,b+h)^2+2E N(b,b+h)N[a,b]\\
&\sim E N(b,b+h)+2E N(b,b+h)N[a,b]\\
&\sim hE\f(V(b))+2hEN[a,b]\f(V(b)-b),\,h\downarrow0,
\end{align*}
and hence:
$$
EN[a,b]^2=EN[a,b]+2E\int_a^b N[a,c]\f(V(c)-c)\,dc.
$$
Moreover, using an obvious time reversal argument, also used in \cite{piet:88}, we get:
\begin{align*}
&E\int_a^b N[a,c]\f(V(c)-c)\,dc=E\int_a^b \left\{\int_a^c \f(-V(c_1)+c_1)\,dc_1\right\}\f(V(c)-c)\,dc\\
&=\int_{a<c_1<c_2<b}E\left\{\f(-V(c_1)+c_1)\f(V(c_2)-c_2)\right\} dc_1dc_2.
\end{align*}
Note:
\begin{align*}
&\int_{a<c_1<c_2<b}E\left\{\f(-V(c_1)+c_1)\f(V(c_2)-c_2)\right\} dc_1dc_2\\
&=\int_{a<c_1<c_2<b}\mbox{\rm cov}\left(\f(-V\left(c_1)+c_1\right),\f\left(V(c_2)-c_2\right)\right) dc_1dc_2\\
&\qquad\qquad\qquad\qquad\qquad+\int_{a<c_1<c_2<b}E\f(-V(c_1)+c_1)E\f(V(c_2)-c_2) \,dc_1dc_2\\
&=\int_{a<c_1<c_2<b}\mbox{cov}\left(\f(-V\left(c_1)+c_1\right),\f\left(V(c_2)-c_2\right)\right) dc_1dc_2\\
&\qquad\qquad\qquad\qquad\qquad\qquad\qquad\qquad+\left(E\f(V(0))\right)^2\int_{a<c_1<c_2<b} \,dc_1dc_2\\
&=\int_{a<c_1<c_2<b}\mbox{cov}\left(\f(-V\left(c_1)+c_1\right),\f\left(V(c_2)-c_2\right)\right) dc_1dc_2
+\tfrac12k_1^2(b-a)^2,
\end{align*}
so we get:
$$
\mbox{var}\bigl(N[a,b]\bigr)
=EN[a,b]+2\int_{a<c_1<c_2<b}\mbox{cov}\left(\f(-V\left(c_1)+c_1\right),\f\left(V(c_2)-c_2\right)\right) dc_1dc_2.
$$
\eop

\vspace{0.3cm}
We prove in the sequel that the dependence between $\f(-V\left(c_1)+c_1\right)$ and $\f\left(V(c_2)-c_2\right)$ dies out exponentially fast, as $c_2-c_1\to\infty$, which, together with part (iii) of Lemma \ref{lemma:p}, gives that the covariance of $\f(-V\left(c_1)+c_1\right)$ and $\f\left(V(c_2)-c_2\right)$ dies out exponentially fast, as $c_2-c_1\to\infty$. Hence we get:
$$
\mbox{var}\bigl(N[a,b]\bigr)\sim k_2(b-a),\,b-a\to\infty,
$$
for a constant $k_2\ge0$.

There are several ways in which one could try to determine the constant $k_2$.
One possible approach is to use an integro-differential equation to determine the constant $k_2$, following a suggestion on p.\ 546 of \cite{piet:85}. Let the function $k(a,t)$ be defined by
\begin{align*}
k(a,t)&=E\left\{\f(V(0))\bigm|V(a)=t\right\},\,a\le0.
\end{align*}
Then we have, for $a<0$,
\begin{align*}
k(a,t)&=E\left\{\f(V(0))\bigm|V(a)=t\right\}\\
&=E\left\{E\left\{\f(V(0))\bigm|V(a+h)\right\}\bigm|V(a)=t\right\}\\
&= h\int_{u=0}^{\infty} k(a,t+u)\frac{2up(u)g(t-a+u)}{g(t-a)}\,du\\
&\qquad\qquad+k(a+h,t)\left\{1-h\int_{u=0}^{\infty}\frac{2up(u)g(t-a+u)}{g(t-a)}\,du\right\}+o(h)\\
&=k(a+h,t)+h\int_{u=0}^{\infty} \left\{k(a,t+u)-k(a,t)\right\}\frac{2up(u)g(t-a+u)}{g(t-a)}\,du+o(h),\,h\downarrow0.
\end{align*}
Hence we get:
$$
\frac{\partial}{\partial a}k(a,t)
=-\int_{u=0}^{\infty} \{k(a,t+u)-k(a,t)\}\frac{2up(u)g(t-a+u)}{g(t-a)}\,du,
$$
which leads to the integral equation
\begin{align*}
k(a,t)=k(0,t)+\int_a^0\,db \int_{u=0}^{\infty}\{k(b,t+u)-k(b,t)\}\frac{2up(u)g(t-b+u)}{g(t-b)}\,du,\,a\le0.
\end{align*}
Note that
$$
E\left\{\f(V(0))\bigm|V(a)-a=t\right\}=k(a,t+a).
$$
Using this approach,
the constant $k_2$ was approximated numerically, 
on a grid with stepsize $10^{-3}$ in both coordinates on the interval $[-10,10]$, using the boundary condition
$$
k(0,t)=\f(t),\,t\in\R,
$$
and replacing integrals by Riemann sums.
In this way we obtained:
$$
2\int_0^{\infty}\mbox{cov}\bigl(\f(V(0)),\f(-V(a)+a)\bigr) da\approx-1.11891,
$$
which would give: $k_2\approx0.986$.
However, since the numerical computations seemed somewhat unstable, we have more trust in the value obtained by simulating the vertex process directly, without first generating Brownian motion, in the way described below.

One could also try to determine an approximate value of the constant $k_2$ by simulating Brownian motion directly. However, since one needs very long intervals (or, alternatively, a rescaling which also leads to very computation-intensive simulation), it is doubtful that we get a good approximation in this way. See also the discussion in Remark \ref{remark:Woodroofe} on the constant $k_1(1/2)$, obtained by simulating Brownian motion directly in \cite{meyer_woodroofe:00}, which gave the value $1.289$, while the analytically determined value is $1.32826$.

We can use Theorem \ref{th:jump_measure} to generate the process $\{V(a):a\in\R\}$ without first generating Brownian motion. This method of generating the vertices was also used in \cite{nagaev:95} and \cite{piet:11}, for generating the vertices of convex hulls of Poisson processes of points in the plane, and seemed to work rather well in that situation.

We start the process at time zero, by generating $V(0)$ according to the ``Chernoffian" distribution $f_{V(0)}$, given by:
$$
f_{V(0)}(x)=\tfrac12g(-x)g(x),\,x\in\R,
$$
where $g$ is defined as in Theorem \ref{th:jump_measure}. Suppose this gives $V(0)=x$. Next we generate the waiting time until a jump according to the distribution function
$$
F_x(a)=1-\exp\left\{-\int_{b=0}^a\f(x-b)\,db\right\}=1-\exp\left\{-\int_{u=x-a}^{x}\f(u)\,du\right\},\,a>0.
$$
where $\f(u)$ is the integrated jump measure, starting from position $u$. Suppose that this gives the jump time $a>0$. Then we generate a jump according to the jump density
$$
u\mapsto \frac{2g(x-a+u)up(u)}{g(x-a)\f(x-a)}\,,\,u>0.
$$
This gives a new position $y$ from which we generate a waiting time, according to the distribution function
$$
F_{y-a}(b-a)=1-\exp\left\{-\int_{u=y-b}^{y-a}\f(u)\,du\right\},\,b-a>0,
$$
which gives a new jump time $b>a$, from which we generate a jump according to the jump density
$$
u\mapsto \frac{2g(y-b+u)up(u)}{g(y-b)\f(y-b)}\,,\,u>0,
$$
and so on. Defining
$$
\Phi(x)=\left\{\begin{matrix}
\int_{u=0}^x\f(u)\,du &,\,x\ge0,\\
&\\
-\int_{u=x}^0\f(u)\,du &,\,x<0,
\end{matrix}
\right.
$$
we can write:
$$
F_{y-a}(b-a)=1-\exp\left\{-\left\{\Phi(y-a)-\Phi(y-b)\right\}\right\}.
$$

We indeed used this procedure to generate the process $V$. Instead of the jumps lengths themselves we generated the square roots of the jump lengths, which have a bounded density, in contrast with the jump lengths, which have a density which is unbounded near zero. This enabled us to generate the square roots of the jump lengths by rejection sampling, since we can compute the density from the theory above (but note that this density depends on the value of $x-a$, so we get a family of densities, parametrized by $x-a$). 

The waiting times between jumps can be generated using the following observations. Note that, for a uniform random variable $U$:
\begin{align*}
&F_x(u)=\P\left\{U\le F_x(u)\right\}=\P\left\{-\log\{1-U\}\le-\log\{1-F_x(u)\}\right\}\\
&=\P\left\{-\log\{1-U\}\le \Phi(x)-\Phi(x-u)\right\}=\P\left\{\Phi_x^{-1}(W)\le u\right\},
\end{align*}
where $W$ is standard exponentially distributed and $\Phi_x^{-1}$ is the inverse of the function
$$
u\mapsto \Phi(x)-\Phi(x-u),\,u\ge0.
$$
Hence the waiting times between jumps can be generated by generating the random variables $\Phi_{x-a}^{-1}(W)$, where $W$ has a standard exponential distribution; $\Phi_{x-a}^{-1}$ was computed on a equidistant grid, with distance $10^{-3}$ between successive points of the grid, and with linear interpolation between points of the grid. In this way we found in $10,000$ simulations, where $a$ ran through the interval $[0,10^4]$: $k_1\approx2.1082$ (note that this is very close to the analytically determined value $2.10484$) and  $k_2\approx1.029$.

The alternative characterization of the jump process, used in the simulations, is given in the following theorem.

\begin{theorem}
\label{th:jump_measure2}
The process $\{V(a):a\in\R\}$ is a Markovian pure jump process, where the jump density at time $a$ is given by
\begin{equation}
\label{jumpdens2}
u\mapsto \frac{2g(x-a+u)up(u)}{g(x-a)\f(x-a)}\,,\,u>0.
\end{equation}
given $V(a-)=x$, and where the distribution function of the waiting time till the next jump is given by
\begin{equation}
F_{x-a}(b-a)=1-\exp\left\{-\int_{u=x-b}^{x-a}\f(u)\,du\right\},\,b-a>0,
\end{equation}
given $V(a)=x$.
\end{theorem}

\begin{remark}
{\rm By part (iii) of Lemma \ref{lemma:p}, we have:
$$
\f(-u)\sim2u^2,\,u\to\infty.
$$
This yields, for fixed $x,a\in\R$,
$$
\int_{u=x-b}^{x-a}\f(u)\,du\sim \tfrac23(b-x)^3\sim \tfrac23b^3,\,b\to\infty,
$$
implying
$$
\log\left\{1-F_{x-a}(b-a)\right\}\sim -\tfrac23b^3,\,b\to\infty.
$$
This is in accordance with:
$$
\log\left(1-\P\left\{|V(a)-a|>t\right\}\right)\sim-\tfrac23t^3,\,t\to\infty
$$
see Corollary 3.4, part (iii), in \cite{piet:89}.
}
\end{remark}

\vspace{0.3cm}
We summarize our findings on the variance in the following lemma.

\begin{lemma}
\label{lemma:variance}
Let $N[a,b]$ be the number of jumps of the process $V$ in the interval $[a,b]$. Then
$$
\mbox{\rm var}\bigl(N[a,b]\bigr)\sim k_2(b-a),\mbox{ as }b-a\to\infty,
$$
where
$$
k_2=k_1+2\int_{-\infty}^0\mbox{\rm cov}\bigl(\f(-V(b)+b),\f(V(0)\bigr)\,db\approx1.029,
$$
and $\f$ is defined by (\ref{phi_b}). The value of the constant $k_2$ was determined by simulating the vertex process directly in the way described above, using Theorem \ref{th:jump_measure2}, by $10^4$ runs on the interval $[0,10^4]$. 
\end{lemma}

For the central limit result, we also need the following lemma.

\begin{lemma}
\label{lemma:mixing}
The process $V(a):a\in\R\}$ is strongly mixing with strong mixing function
$$
\a(d)=c\exp\left\{-\tfrac1{12}d^3\right\},
$$
for a constant $c>0$. More specifically, for arbitrary $a\in\R$ we have:
$$
\sup |\P(A\cap B)-\P(A)\P(B)|\le c \exp\left\{-\tfrac1{12}d^3\right\},
$$
for all $A\in\s\{V(b):b\le a\}$ and $B\in\s\{V(b):b\ge a+d\}$.
\end{lemma}

\noindent
{\bf Proof.} The proof proceeds along similar lines as the proof of Theorem 3.3 in \cite{piet_rik_gerard:99}. Consider, for $a_1,\dots,a_k\le a$ and $b_{\ell}\ge\dots b_1\ge a+d$, the events
$$
E_1=\left\{V(a_1)\in A_1\,\dots,V(a_k)\in A_k\right\},\qquad
E_1=\left\{V(b_1)\in B_1\,\dots,V(b_{\ell})\in B_{\ell}\right\},
$$
for Borel sets $A_1,\dots,A_k$ and $B_1,\dots,B_{\ell}$. Define
$$
M=\tfrac12d,\qquad V^M(b)=\mbox{argmax}_{|t-b|\le M}\left\{W(t)-(t-b)^2\right\},
$$
and consider the events
$$
E_1'=E_1\cap\left\{V(a)=V^M(a)\right\},\qquad E_2'=E_2\cap\left\{V(a+d)=V^M(a+d)\right\}.
$$
By monotonicity, the event $E_1'$ only depends on the increments of Brownian motion before time $a+M$, and the event $E_2'$ only depends on the increments of Brownian motion after time $a+d-M$. By the definition of $M$ and the independent increments property of Brownian motion, this implies that the events $E_1'$ and $E_2'$ are independent, and hence
\begin{equation}
\label{event_independence}
\P\left(E_1'\cap E_2'\right)=\P(E_1')\P(E_2').
\end{equation}
Furthermore, by Corollary 3.4 of \cite{piet:89} we get:
\begin{align*}
&\P\left\{E_1'\ne E_1\right\}\le \P\left\{V(a)\ne V^M(a)\right\}
\le2\P\{V(a)>a+M\}\\
&=2\P\{V(0)>M\}\sim \frac{2^{8/3}}{\ai'(\tilde a_1)}\int_M^{\infty}|t|\exp\left\{-\tfrac23|t|^3+2^{1/3}\tilde a_1|t|\right\}\,dt\\
&\le \frac{2^{8/3}}{\ai'(\tilde a_1)}\int_M^{\infty}|t|\exp\left\{-\tfrac23|t|^3\right\}\,dt
\sim\frac{2^{5/3}}{M\ai'(\tilde a_1)}e^{-\tfrac23M^3}\le 6e^{-\tfrac23M^3},\,M\to\infty.
\end{align*}
where, as before, $\ai$ is the Airy function $\ai$ and $\tilde a_1$ its largest zero on the negative halfline.
The probability $\P\{E_2\ne E_2'\}$ can be handled in a similar way.

Hence we get:
\begin{align*}
\left|\P\left(E_1\cap E_2\right)-\P(E_1)\P(E_2)\right|
\le \P\left\{E_1'\ne E_1\right\}+\P\left\{E_2'\ne E_2\right\}\le ce^{-\tfrac23M^3},
\end{align*}
for a constant $c>0$, and the result follows.\eop

\vspace{0.3cm}
The following lemma shows that all moments of $N[a,b]$ exist, for all $b>a$.

\begin{lemma}
\label{lemma:mgf}
We have:
$$
E e^{\l N[a,b]}<\infty,
$$
for all $\l>0$ and all $b>a$.
\end{lemma}

\noindent
{\bf Proof.}
By the stationarity it is sufficient to prove this for $N[0,a]$. For $\l>0$, the process
$$
a\mapsto\exp\left\{\l N[0,a]-\int_0^a \left\{e^{\l}-1\right\}\f(V(b)-b)\,db\right\},\,a\ge0,
$$
is a martingale w.r.t.\ the filtration $\{{\cal F}_a:a\ge0\}$, where
$$
{\cal F}_a=\s\{V(b),\,b\in[0,a]\},\,a\ge0.
$$
So we find
$$
E\exp\left\{\l N[0,a]-\int_0^a \left\{e^{\l}-1\right\}\f(V(b)-b)\,db+4a\left\{e^{\l}-1\right\}V(0)^2\right\}
=Ee^{4a\left\{e^{\l}-1\right\}V(0)^2},
$$
for each $a\ge0$. Moreover, since, according to part (iii) of Lemma \ref{lemma:p},
$$
\f(t)\sim 2t^2,t\to-\infty,\qquad \f(t)\sim t^{-1},\,t\to\infty,
$$
we have that
\begin{align*}
&4a\left\{e^{\l}-1\right\}V(0)^2-\int_0^a \left\{e^{\l}-1\right\}\f(V(b)-b)\,db\\
&=\int_0^a \left\{e^{\l}-1\right\}\left\{4V(0)^2-\f(V(b)-b)\right\}\,db
\end{align*}
is bounded below, say by $-M$, where $M\ge0$, using $V(0)\le V(b)$, for $b\ge0$.
Hence:
\begin{align*}
&E\exp\left\{\l N[0,a]\right\}
\le E\exp\left\{M+4a\left\{e^{\l}-1\right\}V(0)^2\right\}<\infty,
\end{align*}
for all positive $a$ and $\l$, since, by part (iii) of Corollary 3.4 of \cite{piet:89},
$$
\int e^{\a x^2}f_{V(0)}(x)\,dx<\infty,
$$
for all $\a>0$.\eop

\vspace{0.3cm}
We are now ready to prove our main result.

\vspace{0.3cm}
\noindent
{\bf Proof of Theorem \ref{th:asymp_normality_BM}.} By stationarity, we only have to prove the result for the interval $[0,n]$. We  have:
$$
N[0,n]=\sum_{k=1}^n N[k-1,k],
$$
where the $N[k-1,k],\,k=1,2,\dots$ form a stationary sequence. By Lemma \ref{lemma:mgf}, all moments of $N[0,1]$ exist. This fact, together with the mixing condition of Lemma \ref{lemma:mixing} imply the result, using, e.g., Theorem 18.5.3 of \cite{ibra_lin:71}.
\eop

\section{The jumps of the Grenander estimator}
\label{section:Grenander}
As an application of the results in section \ref{section:BM} we now discuss the use of these results in deriving the asymptotic normality of the number of jumps of the Grenander estimator $\hat f_n$ of a strictly decreasing density $f_0$ on $[0,M]$, $M>0$. Note that the number of jumps of the Grenander estimator is the same as the number of segments of the least concave majorant of the empirical distribution function, which, in turn, is the same as the number of vertices of the least concave majorant minus one.
To keep the length of the present paper within reasonable bounds, we only give a sketch of the proof. Full details will be given elsewhere.

We have the following result.

\begin{lemma}
\label{lemma:Grenander_asymp}
Let $N_n$ the number of jumps of $\hat f_n$, where $\hat f_n$ is the Grenander estimator, based on a sample of size $n$ from $f_0$, where $f_0$ is a decreasing continuous density which stays away from zero on its support $[0,M]$,  with a continuous derivative $f_0'$, which also stays away from zero on $[0,M]$, where one-sided derivatives are taken at the endpoints. Then:
$$
EN_n\sim k_1 n^{1/3}\int_0^M \left|f_0'(x)^2/(4f_0(x))\right|^{1/3}\,dx,\,n\to\infty,
$$
and
$$
\mbox{\rm var}(N_n)\sim k_2 n^{1/3}\int_0^M \left|f_0'(x)^2/(4f_0(x))\right|^{1/3}\,dx,\,n\to\infty,
$$
where the constants $k_1$ and $k_2$ are defined as in Theorem \ref{th:asymp_normality_BM}, that is:
$$
k_1\approx2.10848\qquad \mbox{\rm and }\qquad k_2\approx1.029.
$$
\end{lemma}

\noindent
{\bf Sketch of proof.} Note that, as in \cite{piet:85}, p.\ 542, we can introduce locally a process $V_n$, defined by
$$
V_n(s)=c_2n^{1/3}
\sup\left\{u:\F_n(t+u)-\F_n(t)-(a+n^{-1/3}c_1s)u\mbox{ is maximal}\right\},
$$
where
$$
c_1=2\{\tfrac12a|f_0'(g_0(a))|\bigr\}^{1/3}=\{4a|f_0'(g_0(a))|\bigr\}^{1/3}\qquad\mbox{ and }\qquad c_2=\left\{\frac{f_0'(g_0(a))^2}{4a}\right\}^{1/3},
$$
and $g_0$ is the inverse of $f_0$.
Here $\F_n$ is the empirical df and $a=f_0(t)$ for an interior point $t$ of the support of $f_0$. The process $V_n$ is the (local) inverse of the slope process. As noted in \cite{piet:85}, the process $V_n$ converges in distribution in the Skorohod topology to the process $V$, where $V$ is the process of locations of maxima, discussed in section \ref{section:BM} (where $c_1$ has an extra factor $2$, to obtain $V(s)$ instead of $V(\tfrac12s)$ in the limit).

The jumps of the limiting process are a stationary locally finite point process, implying that the number of jumps of the process $V_n$ on an interval $[b,c]$ converges in distribution to the number of jumps of $V$ on the same interval. Hence, defining
$$
\tilde V_n(s)=V_n(c_1^{-1}s),
$$
we get that the number of jumps of $\tilde V_n$ on an interval $[b,c]$ converges in distribution to the number of jumps of $V$ on an interval of length $c_1^{-1}(c-b)$.

Now note that
$$
\tilde V_n(s)=c_2n^{1/3}
\sup\left\{u:\F_n(t+u)-\F_n(t)-(a+n^{-1/3}s)u\mbox{ is maximal}\right\},
$$
which has, on an interval $[b,c]$, the same number of jumps as the process $U_n$, defined by
$$
U_n(\a)=\sup\left\{x\ge0:\F_n(x)-\a x\mbox{ is maximal}\right\},\,\a>0.
$$
on the interval $[a+bn^{-1/3},a+cn^{-1/3}]$.

We can strengthen this argument somewhat (full details will be given elsewhere), using a Poissonization argument together with a strong approximation result of \cite{kurtz:78}, to show that the expectation and variance of the number of jumps of $U_n$ on an interval
$$
[a-n^{-1/3}\log n,a+n^{-1/3}\log n],
$$
are of order
$$
2k_1c_1^{-1}n^{1/3}\log n\qquad\mbox{ and }\qquad 2k_2c_1^{-1}n^{1/3}\log n,
$$
respectively, where $k_1$ and $k_2$ are defined as in Theorem \ref{th:asymp_normality_BM}.

Partitioning the interval $[f_0(M),f_0(0)]$ into $K_n$ intervals of length of order $2n^{-1/3}\log n$, with midpoints $a_j$, we get:
\begin{align*}
EN_n&\sim k_1n^{1/3}\sum_{j=1}^{K_n}\bigl\{4a_j|f_0'(g_0(a_j))|\bigr\}^{-1/3}\,\,2n^{-1/3}\log n\\
&\sim k_1n^{1/3}\int_{f_0(M)}^{f_0(0)}\left\{4a|f_0'(g_0(a))|\right\}^{-1/3}\,da
=k_1n^{1/3}\int_0^M\left\{4f_0(x)|f_0'(x)|\right\}^{-1/3}\left|f_0'(x)\right|\,dx\\
&=k_1n^{1/3}\int_0^M\left\{\frac{f_0'(x)^2}{4f_0(x)}\right\}^{1/3}\,dx.
\end{align*}
A similar argument, again using Riemann sums approximating the corresponding integral, gives the result for the variance.\eop

\vspace{0.3cm}
The conditions in Lemma \ref{lemma:Grenander_asymp} are probably stronger than needed, and we give two examples below which may also satisfy the result, but do not satisfy the conditions of the lemma. In the first 
example $f_0$ does not stay away from zero on $[0,M]$, and in the second example $f_0$ has infinite support.

If $f_0(x)=2(1-x)$ on $[0,1]$, we get:
$$
\int_0^1 \left|f_0'(x)^2/(4f_0(x))\right|^{1/3}\,dx=3\cdot4^{-2/3}\approx1.19055,
$$
and hence:
$$
EN_n\sim 2.51\,n^{1/3},\qquad
\mbox{var}(N_n)\sim 1.225\,n^{1/3},\,n\to\infty.
$$
A simulation of $1000$ samples with $n=1000$ gave as mean number of the number of jumps $N_{n,k}$, $k=1,\dots,1000$,
$$
\bar N_{n}=\frac{\sum_{k=1}^{1000}N_{n,k}}{1000} = 2.5026\,n^{1/3}.
$$
and as variance
$$
\frac{\sum_{k=1}^{1000}\left(N_{n,k}-\bar N_n\right)^2}{999} = 1.238\,n^{1/3}.
$$

If $f_0$ is the standard exponential density, we get:
$$
\int_0^1 \left|f_0'(x)^2/(4f_0(x))\right|^{1/3}\,dx=3\cdot2^{-2/3}\approx1.88988,
$$
and hence:
$$
EN_n\sim 3.98477\,n^{1/3},\qquad
\mbox{var}(N_n)\sim 1.945\,n^{1/3},\,n\to\infty,
$$
whereas $1000$ samples of size $n=1000$ yielded:
$$
\bar N_{n}=\frac{\sum_{k=1}^{1000}N_{n,k}}{1000} = 3.64544\,n^{1/3}.
$$
and as variance
$$
\frac{\sum_{k=1}^{1000}\left(N_{n,k}-\bar N_n\right)^2}{999} = 1.86570\,n^{1/3}.
$$

\vspace{0.3cm}
Analogously to the result for Brownian motion, we have:

\begin{theorem}
Let $N_n$ the number of jumps of $\hat f_n$, where $\hat f_n$ is the Grenander estimator, based on a sample of size $n$ from $f_0$, where $f_0$ satisfies the conditions of Lemma \ref{lemma:Grenander_asymp}. Moreover, let the constants $k_1$ and $k_2$ be defined as in Theorem \ref{th:asymp_normality_BM}. Then:
$$
n^{-1/6}\left\{N_n-k_1 n^{1/3}\int_0^M \left\{f_0'(x)^2/(4f_0(x))\right|^{1/3}\,dx\right\}\stackrel{\cal D}\longrightarrow N(0,\s^2),
$$
where $N(0,\s^2)$ is a normal distribution with expectation zero and variance
$$
\s^2=k_2\int_0^M \left|f_0'(x)^2/(4f_0(x))\right|^{1/3}\,dx.
$$
\end{theorem}

Although the proof proceeds along similar lines as the proof of Theorem \ref{th:asymp_normality_BM} in section \ref{section:BM}, the embedding into Brownian motion needs some careful attention, and therefore the details of the proof will be given elsewhere.

\section{Concluding remarks}
\label{section:conclusion}
In the preceding, a central limit result was proved for the number of vertices in an increasing interval of the concave majorant of the process $\{W(t)-t^2,\,t\in\R\}$, where $W$ is two-sided standard Brownian motion, originating from zero. The central limit result involves two constants $k_1$ and $k_2$ for the mean and variance, respectively, see Theorem \ref{th:asymp_normality_BM}. The constant $k_1$ has several representations, for example
$$
k_1=8EV(0)^2=\tfrac83E\max_{t\in\R}\left(W(t)-t^2\right)\approx2.10848,
$$
see (\ref{2nd_representation_k1}), where $a\mapsto V(a)$ is the process of locations of maxima of $W(t)-(t-a)^2$, as a function of $a$. From \cite{meyer_woodroofe:00} we get:
$$
k_1=2E\tilde X(0)+E\tilde X'(0)^2=2E\tilde X(0)+4EV(0)^2
$$
where $\tilde X$ is the concave majorant of the process $\{W(t)-t^2,\,t\in\R\}$. This implies as a side result the relation
$$
E\tilde X(0)=\tfrac23E\max_{t\in\R}\left\{W(t)-t^2\right\},
$$
see Remark \ref{remark:Woodroofe2}.

Much less is known about the constant $k_2$. We used a direct simulation of the vertex process to determine this constant, but there is  room for improvement here. The approximate value we found is close to $1$, and our preliminary value is: $k_2=1.029$.
The basis for the simulation of the vertex process directly, without first generating Brownian motion, is given in Theorem \ref{th:jump_measure2}, which gives the distribution of the (non-exponential) waiting times between jumps of the process $V$, as a function of $a$ and $x$, where $V(a)=x$, and the density of the size of the jumps. The square roots of the jump lengths were generated by rejection sampling, and the waiting times between jumps by generating standard exponential random variables, and by applying the inverse of the cumulative hazard function of the waiting time distribution (again parametrized by $x$ and $a$) on these.

A similar technique of generating vertices of a convex hull was used in \cite{nagaev:95} and \cite{piet:11}, where convex hulls of random points in the plane were studied. The behavior of the least concave majorant of Brownian motion minus a parabola has some remarkable analogies with the behavior of the convex hulls of points chosen uniformly from the interior of a circle, where we also get central limit theorems for the number of vertices, with an expectation and variance which are also both of order $n^{1/3}$, if $n$ is the number of points chosen. On the other hand, the behavior of the concave majorant of one-sided Brownian motion without drift has analogies with the behavior of the convex hulls of points drawn uniformly from a convex polygon, where we  get central limit theorems for the number of vertices with an asymptotic expectation and variance which are both of order $\log n$, if $n$ is the number of points chosen, as shown in \cite{piet:88}.

\vspace{0.3cm}
\noindent
{\bf Acknowledgements}
I am grateful to the referees  for their careful reading and useful comments.

\bibliographystyle{amsplain}
\bibliography{references}

\end{document}